\documentclass[12pt]{article}

\usepackage{amssymb,amsmath,latexsym,amsthm}
\usepackage[english]{babel}

\newtheorem{Theorem}{Theorem}[section]
\newtheorem{lemma}{Lemma}[section]
\newtheorem{proposition}{Proposition}[section]
\newtheorem{cor}{Corollary}[section]

\theoremstyle{definition}

\newtheorem*{example}{Example}

\begin{document}

\title{Maximal selectivity for orders in fields}


\author{\sc Luis Arenas-Carmona}


\newcommand\Q{\mathbb Q}
\newcommand\alge{\mathfrak{A}}
\newcommand\Da{\mathfrak{D}}
\newcommand\Ea{\mathfrak{E}}
\newcommand\Ha{\mathfrak{H}}
\newcommand\oink{\mathcal O}
\newcommand\matrici{\mathbb{M}}
\newcommand\Txi{\lceil}
\newcommand\ad{\mathbb{A}}
\newcommand\enteri{\mathbb Z}
\newcommand\finitum{\mathbb{F}}
\newcommand\bbmatrix[4]{\left(\begin{array}{cc}#1&#2\\#3&#4\end{array}\right)}

\maketitle

\begin{abstract}
If $\Ha\subseteq\Da$ are two orders in a central simple algebra
$\alge$ with $\Da$ of maximal rank, the theory of representation
fields describes the set of spinor genera of orders in the genus
of $\Da$ representing the order $\Ha$. When $\Ha$ is contained in
a maximal subfield of $\alge$ and the dimension of $\alge$ is the
square of a prime $p$, the proportion of spinor genera
representing $\Ha$ has the form $r/p$, in fact, when the
representation field exists, this proportion is either $1$ or
$1/p$. In the later case the order $\Ha$ is said to be selective
for the genus. The condition for selectivity is known when $\Da$
is maximal and also when $p=2$ and $\Da$ is an Eichler order. In
this work we describe the orders $\Ha$ that are selective for at
least one genus of orders of maximal rank in $\alge$.
\end{abstract}

\bigskip
\section{Introduction}

Let $K$ be a number field. Let $\alge$ be a central simple algebra (CSA) over $K$,
and let $\Da$ be an order of maximal rank in $\alge$. Let $\Sigma$ be the spinor
class field for the genus $\mathbb{O}$ of $\Da$ as defined in \cite{spinor}. In
particular, $\Sigma/K$ is an abelian extension such that the spinor genera in
$\mathbb{O}$ can be described in terms of the Galois group
$\mathcal{G}=\mathrm{Gal}(\Sigma/K)$. More precisely, there exists a map
$$\rho:\mathbb{O}\times\mathbb{O}\rightarrow\mathcal{G},$$
Such that $\Da'$ belongs to the spinor genus $\mathrm{spin}(\Da)$
if and only if $\rho(\Da,\Da')=\mathrm{Id}_\Sigma$. Furthermore,
the map $\rho$ satisfies
$\rho(\Da,\Da'')=\rho(\Da,\Da')\rho(\Da',\Da'')$ for any triple
$(\Da,\Da',\Da'')\in\mathbb{O}^3$ (\cite{spinor}, \S3). Note that
when strong approximation applies to the algebraic group
$\mathrm{SL}_1(\alge)$, spinor genera coincide with conjugacy
classes \cite{spinor}. All of the above generalize to arbitrary
global fields. This is done in (\cite{abelianos}, \S2) when $\Da$
is maximal, although the last hypothesis is not actually used
there. In fact, everything that follows holds in the function
field case as much as in the number field case. It also holds for
$S$-orders for an arbitrary set of places $S$ containing the
archimedean places if any.

For any suborder $\Ha$ of $\Da$, we can define two intermediate
subfields:
\begin{enumerate}\item
The upper representation field $F=F^-(\Da|\Ha)$ for the pair
$(\Da|\Ha)$ is the smallest subfield of $\Sigma$ containing $K$
such that the order $\Ha$ is represented by the spinor genus
$\mathrm{Spin}(\Da')$ for every order $\Da'\in\mathbb{O}$
satisfying $\rho(\Da,\Da')|_F=\mathrm{Id}_F$. \item The lower
representation field $F=F_-(\Da|\Ha)$ is the largest subfield of
$\Sigma$ such that $\Ha\subseteq\Da'$ implies
$\rho(\Da,\Da')|_F=\mathrm{Id}_F$.
\end{enumerate}
When $F_-(\Da|\Ha)=F^-(\Da|\Ha)$, we call this field the
representation field $F(\Da|\Ha)$ for the pair $(\Da|\Ha)$. In
this case we say that the representation field exists. The
existence of a representation field $F$ for $\Ha$ implies that the
proportion of conjugacy classes in $\mathbb{O}$ representing $\Ha$
is $[F:K]^{-1}$. This fact was first studied by Chevalley
\cite{Chevalley} when $\alge$ is a matrix algebra of arbitrary
dimension, $\Da$ is a maximal order, and $\Ha$ is the maximal
order in a maximal subfield of $\alge$. Later computations for
$F(\Da|\Ha)$ can be summarized in the following table:

\footnotesize
\[
\begin{tabular}{ | c | c | c | c | c | c |c| }
  \hline
  Year & Ref. & $\alge$ (CSA) & $\Da$ & $\Sigma/K$ & $\Ha$ & $K\Ha$ \\
  \hline\hline
  1936 & \cite{Chevalley} & matrix & max. & & max. & field (max.) \\  \hline
  1999 & \cite{FriedmannQ} & quaternion & max. & &  & commutative  \\  \hline
  2003 & \cite{spinor} & NPR & max. & & max. & field (maximal) \\  \hline
  2004 & \cite{Guo},\cite{Chan} & quaternion & EO & & & commutative  \\   \hline
    2004 & $\cite{Chan}^{e}$ & quaternion & & & & commutative  \\   \hline
  2008 & $\cite{eichler}^{e}$  &  & max.& $\mathcal{G}^2=\{\mathrm{id}\}$ &  &   \\
  \hline
   2008 & $\cite{eichler}^{c}$  &  & max.& $\mathcal{G}^2\neq\{\mathrm{id}\}$ & GEO & $=\alge$   \\
  \hline
  2010 & \cite{lino2} & prime degree & max. & &  & field  \\  \hline
   2011 & \cite{lino1} & quaternion & & unramified  & & commutative  \\ \hline
  2011 & \cite{abelianos} &  & max. & &  & commutative  \\  \hline
    2011 & \cite{cyclic} &  & max. & & LCMO & CSA  \\
    \hline
\end{tabular}
\]
\normalsize
 NPR above stands for \emph{no partial ramification}, a
weaker condition than prime degree, while LCMO means \emph{locally
cyclic or maximal order}, and (G)EO means \emph{(generalized)
Eichler order}. Here $[n]^c$ denotes a counter-example, i.e., an
example where the representation field fails to exist, while
$[n]^e$ denotes an existencial result. The proof of the main
result in \cite{abelianos} seems to be easy to generalize to other
families of suborders. For example, it is very simple now to write
a general formula for the representation field when $K\Ha$ is
contained in a quaternion algebra. However, the condition that
$\Da$ is maximal is essential in this computations and a
generalization to arbitrary orders $\Da$ of maximal rank seems
unlikely at this point.

In this paper we focus on the case where $\Ha$ is an order of
maximal rank in a maximal subfield of $\alge$. Instead of trying
to give a general formula for all representation fields
$F(\Da|\Ha)$, we focus on the maximal possible representation
field $F_M(\alge|\Ha)=\max_{\Da\subseteq\alge}F(\Da|\Ha)$, where
$\Da$ runs over the set of all orders of maximal rank for which
$F(\Da|\Ha)$ is defined. It follows easily from formula (7) in
\S4.2 of \cite{spinor} that $F^-(\Da|\Ha)\subseteq L$ for every
order $\Da$ of maximal rank, and therefore also
$F_M(\alge|\Ha)\subseteq L$, as long as this maximum exists. Here
we give a formula for $F_M(\alge|\Ha)$ whenever
$\sqrt{\dim_K(\alge)}=p$ is a prime. In particular we prove the
existence of $F_M(\alge|\Ha)$ in this case. For extensions $L/K$
of prime degree, an order of maximal rank in $L$ is selective (in
the sense defined in \cite{FriedmannQ} on \cite{lino2}) for some
genus of orders of maximal rank in $\alge$ if and only if
$F_M(\alge|\Ha)=L$ (Prop. 3.1). Since the even-dimensional and
odd-dimensional cases are different, we state our results in the
next two theorems:

\begin{Theorem}\label{th1}Let $\alge$ be a quaternion algebra and $L\cong K(\sqrt d)$ a maximal subfield.
Then, for any order $\Ha$ of maximal rank in $L$, we have
$F_M(\alge|\Ha)=L$ if and only if $\alge$ ramifies at exactly the
same set of finite primes as the algebra
$\left(\frac{-1,d}{K}\right)$. Otherwise we have
$F_M(\alge|\Ha)=K$.
\end{Theorem}

When $L/K$ is a Galois extension, we say that an order
$\Ha\subseteq L$ is asymmetrical at a non-split place $\wp$ if
$\sigma(\Ha_\wp)\neq\Ha_\wp$ for some element
$\sigma\in\mathrm{Gal}(L/K)$.

\begin{Theorem}\label{th2} Assume that $\dim_K(\alge)=p^2$ where $p$ is an
odd prime, and let $L\subseteq\alge$ be a maximal subfield. Then,
for any order $\Ha$ of maximal rank in $L$, we have
$F_M(\alge|\Ha)=L$ if and only if $L/K$ is Galois and the order
$\Ha$ is asymmetrical at every finite place that is ramified for
$\alge$. If this condition is not satisfied we have
$F_M(\alge|\Ha)=K$.
\end{Theorem}

\section{A continuity principle}

In all of this section, $K$ is a local field and $\alge$ is a
central simple $K$-algebra. We denote by $x\mapsto|x|$ the
absolute value on $\alge$ or $K$. Note that $\alge$ is locally
compact since it is a finite dimensional vector space over the
locally compact field $K$.

\begin{lemma}
Let $\alge$ and $K$ be as above. Then the conjugation stabilizer
of a maximal order is compact in $\alge^*/K^*$.
\end{lemma}

\begin{proof}
Assume first that $\alge$ be a division algebra. We claim that
$\alge^*/K^*$ is compact. The result follows from the claim since
a division algebra has a unique maximal order (\cite{weil}, Ch. 1,
Thm. 6). Since $e=\Big[|\alge^*|:|K^*|\Big]$ is finite, it
suffices to observe that the kernel of the absolute value is
$N=\overline{B(0;1)}-B(0;1)$, where $B(0;1)$ is the open ball in
$\alge$ centered at $0$,  which is a compact set.

 Assume now that $\alge\cong\matrici_m(\alge_0)$
for some division algebra $\alge_0$. The conjugation-stabilizer of
a maximal order $\Da$ is $\Da^*\alge_0^*/K^*$, whence the
conclusion follows since $\alge_0^*/K^*$ is compact by the
previous result and $\Da^*$ is compact since it is closed in the
compact set $\Da$.
\end{proof}

In any metric space $(X,d)$ we define for every pair of subsets
$A$ and $B$ of $X$
$$\rho(A,B)=\sup_{a\in A}d(a,B).$$
When $B$ is closed, we have $\rho(A,B)=0$ if and only if
$A\subseteq B$. Note that $\rho$ is not a metric, since it is not
symmetric, but $\hat{\rho}(A,B)=\rho(A,B)+\rho(B,A)$ is a metric
on the collection of compact subsets of $X$ called the Hausdorff
metric.

In all that follows, for every pair of orders $\Da$ and $\Ha$ in
$\alge$ we denote $$H(\Da|\Ha)=\{N(u)|u\in\alge^*, u^{-1}\Ha
u\subseteq \Da\},\qquad H(\Ha)=H(\Ha|\Ha),$$ where
$N:\alge^*\rightarrow K^*$ is the reduced norm.

\begin{proposition}
Let $\alge$ be a central simple algebra over the local field $K$.
Assume the order $\Ha$ is contained in finitely many maximal
orders, and let $\{\Da_t\}_{t\in \mathbb{N}}$ be a sequence of
orders such that
$\rho(\Da_t,\Da)\stackrel{t\rightarrow\infty}\longrightarrow0$.
Then, in the set theoretical sense:
$$\mathop{\mathrm{lim}\,\mathrm{sup}}_{t\rightarrow\infty} H(\Da_t|\Ha)\subseteq
H(\Da|\Ha).$$
\end{proposition}

\begin{proof}
It suffices to prove that if $a\in H(\Da_t|\Ha)$ for infinitely
many values of $t$, then $a\in H(\Da|\Ha)$. The hypotheses imply
$a\in N(y_t)K^{*2}$ for some $y_t$ satisfying $y_t\Ha
y_t^{-1}\subseteq\Da_t$. Let $\Da'$ be a maximal order containing
$\Da$. Since $\Da'$ is open, then $\Da_t\subseteq\Da'$ for $t$
sufficiently large. If $y_t\Ha y_t^{-1}\subseteq\Da_t$, then
$\Ha\subseteq y_t^{-1}\Da_t y_t$, whence $\Ha\subseteq
y_t^{-1}\Da' y_t$ for $t$ sufficiently large. It follows that the
set of maximal orders $\{y_t^{-1}\Da' y_t\}_t$ is finite. Write
$\bar{x}$ for the class in $\alge^*/K^*$ of an element
$x\in\alge^*$. As the stabilizer of $\Da'$ in $\alge^*/K^*$ is
compact, the sequence $\{\bar y_t\}_{t\in \mathbb{N}}$ is
contained in a compact set, whence, by taking a subsequence if
needed, we can assume it is convergent in $\alge^*/K^*$ to an
element $\bar y$. In particular, $y\Ha y^{-1}\subseteq\Da$, and
$N(y_t)\in N(y)K^{*2}$ for $t$ sufficiently large. We conclude
that $a\in N(y)K^{*2}$, and the result follows.
\end{proof}

\begin{cor}\
Let $\Ha$ be an order contained in finitely many maximal orders.
Then there exists $\epsilon=\epsilon(\Ha)$ such that whenever
$\Ha\subseteq\Da$ with $\rho(\Da,\Ha)\leq\epsilon$, then
$H(\Da|\Ha)=H(\Ha)$.
\end{cor}

\begin{proof} Note that the set of quadratic clases is finite, so one inclusion follows from the previous lemma. The opposite inclusion
is immediate, since $H(\Da|\Ha)=H(\Da|\Ha)H(\Ha)$  by the general
theory \cite{spinor}.\end{proof}

\begin{proposition}
Assume $\alge$ is a matrix algebra. If $L=K\Ha$ is a maximal
subfield of $\alge$, then $\Ha$ is contained in finitely many
maximal orders of $\alge$.
\end{proposition}

\begin{proof} Since every pair of embeddings of $L$ into $\alge$ are conjugate,
we can identify $L$ with its natural image in
$\mathrm{Aut}_K(L)\cong\mathrm{Aut}_K(K^n)\cong\alge$. The
$\Ha$-invariant lattices in $K^n$ correspond to fractional
$\Ha$-ideals in $L$. It suffices, therefore, to prove that $K^*$
acts on the set of fractional $\Ha$-ideals with finitely many
orbits. Let $\Lambda$ be a fractional ideal. Multiplying by an
element of $K^*$ if needed we can assume that $\Lambda\subseteq
\oink_L$, but $\Lambda$ is not contained in $\pi_K\oink_L$ for a
uniformizing parameter $\pi_K$ of $K$, i.e., there exist some
element $u\in\Lambda\backslash\pi_K\oink_L$. Since $\oink_L$ is a
valuation ring, we have $\pi_K\oink_L\subseteq u\oink_L$. Since
$\Ha$ has maximal rank in $L$, $\pi^N_K\oink_L\subseteq\Ha$ for
some $N$, whence
$\pi^{N+1}_K\oink_L\subseteq\pi^N_K\oink_Lu\subseteq\Ha\Lambda=\Lambda$.
It follows that $\pi^{N+1}_K\oink_L\subseteq\Lambda\subseteq
\oink_L$ and the result follows.
\end{proof}

\begin{cor}\label{coruse}
Let $\Ha$ be an order of maximal rank in a maximal subfield of
$\alge$. Then there exists $\epsilon=\epsilon(\Ha)$ such that
whenever $\Ha\subseteq\Da$ with $\rho(\Da,\Ha)\leq\epsilon$, then
$H(\Da|\Ha)=H(\Ha)$.\qed
\end{cor}

\begin{example}
Assume $\alge$ is a split quaternion algebra and $\Ha$ is an order
in a maximal unramified subfield $L$. It is proved in \cite{Guo}
that $H(\Da|\Ha)=K^*$ whenever $\Da$ is an Eichler order
representing an order in $L$ strictly containing $\Ha$. This
result does not generalizes to arbitrary orders of maximal rank.
For example if
$\Da_k\stackrel{k\rightarrow\infty}\longrightarrow\oink_L$ in the
Hausdorff metric, then for $k$ big enough we have
$H(\Da_k|\Ha)\subseteq H(\oink_L|\Ha)=H(\oink_L)=\oink_K^*K^{*2}$
according to the computations in \S3.
\end{example}

\section{Computation of $F_M$}

In all that follows, $K$ is a global field, $L/K$ is a field extension of prime
degree $p$ and $\alge$ is a $p^2$-dimensional central simple algebra containing a
copy of $L$. We let $J_K$ be the idele group of $K^*$, $\alge_{\ad}$ the adelization
of the algebra $\alge$, and $N:\alge^*_{\ad}\rightarrow J_K$ the adelic reduced
norm. For every pair of global orders $\Ha\subseteq\Da$ we define $H(\Da|\Ha)$ as
the set
$$\{N(u)|u\in\alge_{\ad}^*, u^{-1}\Ha u\subseteq \Da\}=
\left[\prod_{\wp\notin S}H(\Da_\wp|\Ha_\wp)\times\prod_{\wp\in
S}N(\alge_\wp^*)\right]\cap J_K,$$ where $H(\Da_\wp|\Ha_\wp)$ is defined as in the
preceding section, and $H(\Ha)=H(\Ha|\Ha)$. Define the abelian extension
$F_0(\alge|\Ha)$ of $K$ as the class field corresponding to $K^*H(\Ha)$. Note that
$H(\Ha)$ is a group since it is the image under the reduced norm of the
conjugation-stabilizer of $\Ha$. It is immediate from the general theory that
$H(\Ha)H(\Da|\Ha)=H(\Da|\Ha)$ for any order $\Da$ of maximal rank, whence in
particular $F^-(\Da|\Ha)\subseteq F_0(\alge|\Ha)$. It follows that an order $\Ha$
such that $F_0(\alge|\Ha)=K$ cannot be selective for any genus.

Let $\wp$ be a non-split place for $L/K$. Note that the local conjugation-stabilizer
$N_\wp$ of $\Ha_\wp$ fits into a short exact sequence $L_\wp^*\hookrightarrow N_\wp
\twoheadrightarrow \Gamma_\wp$, where $\Gamma_\wp$ is contained in the Galois group
$\mathrm{Gal}(L/K)$. Furthermore, by Skolem-Noether's Theorem, $\Gamma_\wp$ is
trivial only in the following cases:
\begin{enumerate}
\item $L_\wp/K_\wp$ is not Galois.\item $\Ha_\wp$ is asymmetrical.
\end{enumerate}
In any other case $\Gamma_\wp$ is a cyclic group of order $p$.
Note that in particular $H(\Ha)\supseteq N_{L/K}(J_L)$, whence it
follows that $F_0(\alge|\Ha)\subseteq L$.

\begin{proposition}\label{corun}
Let $\alge$ be a central simple $p^2$-dimensional $K$-algebra, where $p$ is a prime.
If $L=K\Ha$ is a maximal subfield of $\alge$, then there exists an order of maximal
rank $\Da$ in $\alge$ satisfying $F(\Da|\Ha)=F_0(\alge|\Ha)$. In particular,
$F_M(\alge|\Ha)=F_0(\alge|\Ha)$ and $\Ha$ is selective for some genus of maximal
orders of maximal rank in $\alge$ if and only if $F_0(\alge|\Ha)\neq K$.
\end{proposition}

\begin{proof}If $L/K$ is not Galois,
the contention $F^-(\Da|\Ha)\subseteq L$ shows that $F(\Da|\Ha)$ is defined and
equals $K$ for any order $\Da$ of maximal rank, and $F_0(\alge|\Ha)=K$ for the same
reason, whence we can assume that $L/K$ is Galois. Let $T$ be the set of all finite
places $\wp$ satisfying one of the following conditions:
\begin{enumerate}
\item $\alge_\wp$ is ramified. \item $L/K$ is inert at $\wp$ and
$\Ha_\wp$ is not maximal in $L_\wp$. \item $L/K$ is ramified at
$\wp$.
\end{enumerate}
For any $\wp\notin T$ we choose $\Da_\wp$ maximal. Let $\mathbb{H}_\wp$ be the
residuual algebra defined in \cite{abelianos}. When $L/K$ splits at $\wp$, every
representation of the residual algebra $\mathbb{H}_\wp$ has dimension $1$, so we
have $H_\wp(\Da|\Ha)=K_\wp^*$ (Lemma 3.4 in \cite{abelianos}). When $L/K$ is inert
at $\wp$, $\alge_\wp$ is unramified, and $\Ha_\wp$ is maximal in $L_\wp$, every
representation of the residual algebra $\mathbb{H}_\wp$ has dimension $p$ and
therefore $H_\wp(\Da|\Ha)=\oink_\wp^*K^{*p}$ (Lemma 3.4 in \cite{abelianos}). In any
case $$H_\wp(\Ha)\subseteq H_\wp(\Da|\Ha)=N_{L_\wp/K_\wp}(L_\wp^*)\subseteq
H_\wp(\Ha)$$ at all finite places $\wp\notin T$. For the places $\wp\in T$, we
choose $\Da_\wp$ satisfying the conclusion of Corollary \ref{coruse}. The result
follows. The last statement is a consequence of this and the discussion at the
beginning of the section.
\end{proof}

\subparagraph{Proof of Theorem \ref{th1}} Note that every order in
a quadratic extension is symmetric.
 Assume $L\subseteq \alge$ and let $c\in\alge$ be a pure
 quaternion satisfying $cac^{-1}=\bar{a}$ for every $a\in L$. Then
$r=c^2$ is a reduced norm from $L_\wp$ if and only if $\alge_\wp$
splits. Assume first that $\alge_\wp$ is a matrix algebra. Then
the reduced norm $N(c)=-r$ is in $N_{L_\wp/K_\wp}(L_\wp^*)$ if and
only if $-1\in N_{L_\wp/K_\wp}(L_\wp^*)$. Assume next that
$\alge_\wp$ is a division algebra. Then, $F_\wp$ is a field,
whence $N_{L_\wp/K_\wp}(L_\wp^*)$ is a subgroup of index $2$ in
$K_\wp^*$. It follows that $-r\in N_{L_\wp/K_\wp}(L_\wp^*)$ if and
only if $-1\notin N_{L_\wp/K_\wp}(L_\wp^*)$. The result follows in
either case. \qed

\subparagraph{Proof of Theorem \ref{th2}} As in the proof of Proposition \ref{corun}
we can assume that $L/K$ is Galois. Let $\sigma$ be a generator of the Galois group
$\mathrm{Gal}(L/K)$ and assume $L\subseteq \alge$. Fix a local place $\wp$ and let
$c\in\alge_\wp$ be an element satisfying $cac^{-1}=\sigma(a)$ for every $a\in
L_\wp$. Such a $c$ exists by Skolem-Noether's Theorem. For any generator $v$ of
$L/K$, we have $c^rv=\sigma^r(v)c^r$, whence the powers of $c$ are eigenvectors of
the map $x\mapsto xv$, and therefore linearly independent over $L$. It follows that
$L$ and $c$ generate $\alge$ and therefore $r=c^p\in K$, since it is central. Note
that $r$ is a reduced norm from $L_\wp$ if and only if $\alge_\wp$ splits
(\cite{invol}, Prop 30.6). The result follows now since $N(c)=r$. \qed


%


\begin{thebibliography}{xx}



\bibitem{spinor}
{\sc L.E. Arenas-Carmona}, \textit{Applications of spinor class
fields: embeddings of orders and quaternionic lattices}, Ann.
Inst. Fourier \textbf{53} (2003), 2021--2038.

\bibitem{eichler}
 {\scshape Luis Arenas-Carmona}. Relative spinor class
fields: A counterexample, \textit{Archiv. Math.} \textbf{91}
(2008), 486-491.


\bibitem{abelianos}
{\sc L.E. Arenas-Carmona}, \textit{Representation fields for
commutative orders}, to appear in Ann. Inst. Fourier.
arXiv:1104.1809v1 [math.NT].

\bibitem{cyclic}
{\sc L.E. Arenas-Carmona}, \textit{Representation fields for
cyclic orders}, Submitted.

\bibitem{Chan}
{\sc W.K. Chan} and {\sc F. Xu}, \textit{On representations of
spinor genera}, Compositio Math. \textbf{140.2} (2004), 287-300.

\bibitem{Chevalley}
{\sc C. Chevalley}, \textit{L'arithm\'etique sur les alg\`ebres de
matrices}, Herman, Paris, 1936.


\bibitem{FriedmannQ}
{\sc T. Chinburg} and {\sc E. Friedman}, \textit{An embedding
theorem for quaternion algebras}, J. London Math. Soc.
\textbf{60.2} (1999), 33-44.


\bibitem{Guo}
{\sc X. Guo} and {\sc H. Qin}, \textit{An embedding theorem for
Eichler orders}, J. Number Theory \textbf{107.2} (2004), 207-214.

\bibitem{invol}
{\sc M.-A. Knus, A. Merkurjev, M. Rost, and J.-P. Tignol}.
\textit{The book of involutions}, AMS Colloquium Pub. \textbf{44},
1998.

\bibitem{lino1}
{\sc B. Linowitz}, \textit{Selectivity in quaternion algebras},
Preprint.

\bibitem{lino2}
{\scshape B. Linowitz and T.R. Shemanske}, \textit{Embedding
orders into central simple algebras}, To appear in J. Th\'eor.
Nombres Bordeaux.



\bibitem{weil}
{\sc A. Weil}. \textit{Basic Number Theory}, $2^{\mathrm{nd}}$
Ed.,
 Springer Verlag, Berlin, 1973.

\end{thebibliography}
\end{document}